\documentclass[letterpaper,12pt,reqno,twoside]{amsart}

\usepackage{amssymb}
\usepackage{amsmath}
\usepackage[dvips]{graphicx}
\usepackage{psfrag,color}
\usepackage[colorlinks]{hyperref}

\newcommand{\be}{\begin{eqnarray}}
\newcommand{\ee}{\end{eqnarray}}

\newcommand{\beq}{\begin{equation}}
\newcommand{\eeq}{\end{equation}}
\newcommand{\beqn}{\begin{equation*}}
\newcommand{\eeqn}{\end{equation*}}

\newcommand{\N}{\mathbb{N}}
\newcommand{\defas}{\mathrel{\raise.095ex\hbox{$:$}\mkern-4.2mu=}}
\newcommand{\defasr}{\mathrel{=\mkern-4.2mu\raise.095ex\hbox{$:$}}}

\newcommand{\ie}{\textit{i.e.}}

\DeclareMathAlphabet{\mathfat}{U}{bbold}{m}{n}          % Identity operator; requires amsfonts
\newcommand{\one}{\mathfat{1}}

\newtheorem{theo}{Theorem}

\newtheorem{example}[theo]{Example}

\newcommand\cN{{\mathcal N}}

\newcommand\bP{{\mathbb P}}

\newcommand\bR{{\mathbb R}}

\newcommand{\bfP}{{\bf P}}
\newcommand{\bfE}{{\bf E}}

\newcommand{\Var}{\operatorname{Var}}

\begin{document}

\title[Deterministic Walks in Quenched Random Environments...]{Deterministic Walks in Quenched Random Environments of Chaotic Maps}

\author[Tapio Simula]{Tapio Simula}
\email{tapio.simula@gmail.com}
\address[Tapio Simula]{Mathematical Physics Laboratory, Department of Physics, Okayama University, Okayama 700-8530, Japan.}
%\urladdr{}
\author[Mikko Stenlund]{Mikko Stenlund}
\email{mikko@cims.nyu.edu}
\address[Mikko Stenlund]{
Courant Institute of Mathematical Sciences\\
New York, NY 10012, USA; Department of Mathematics and Statistics, P.O. Box 68, Fin-00014 University of Helsinki, Finland.}
\urladdr{http://www.math.helsinki.fi/mathphys/mikko.html}


\date{\today}

\begin{abstract}
This paper concerns the propagation of particles through a quenched random medium. In the  one- and two-dimensional models considered, the local dynamics is given by expanding circle maps and hyperbolic toral automorphisms, respectively. The particle motion in both models is chaotic and found to fluctuate about a linear drift. In the proper scaling limit, the cumulative distribution function of the fluctuations converges to a Gaussian one with system dependent variance while the density function shows no convergence to any function. We have verified our analytical results using extreme precision numerical computations.
\end{abstract}

\maketitle

%%%%%%%%%%%%%%%%%%%%%%%%%%%%
%%%%%%%%%%%%%%%%%%%%%%%%%%%%

\section{Introduction}
Variants of a mechanical model now widely known as the Lorentz gas have occupied the minds of scientists for more than a century. Initially proposed by Lorentz \cite{Lorentz} in 1905 to describe the motion of an electron in a metallic crystal, the model consists of fixed, dispersing, scatterers in $\mathbb{R}^d$ and a free point particle that bounces elastically off the scatterers upon collisions.

If the lattice of scatterers is periodic, the model is also referred to as Sinai Billiards after Sinai, who proved \cite{Sinai70} that the system (with $d=2$) is ergodic if the free path of the particle is bounded. In the latter case it was also proved that, in a suitable scaling limit, the motion of the particle is Brownian \cite{BunimovichSinai-BM}. Sinai's work can be considered the first rigorous proof of Boltzmann's Ergodic Hypothesis in a system that resembles a real-world physical system.

Despite its seeming innocence, the Lorentz gas exhibits a great deal of complexity.  One example is the lack of smoothness of the dynamics caused by tangential collisions of the particle with the scatterers. Another one, the presence of recollisions, is a source of serious statistical difficulties that have not been overcome in the study of the aperiodic Lorentz gas. For more background, see \cite{Tabachnikov,Szasz,ChernovMarkarian,ChernovDolgopyat} and the references therein.

\begin{figure}
   \centering
   \includegraphics[width=0.7\linewidth]{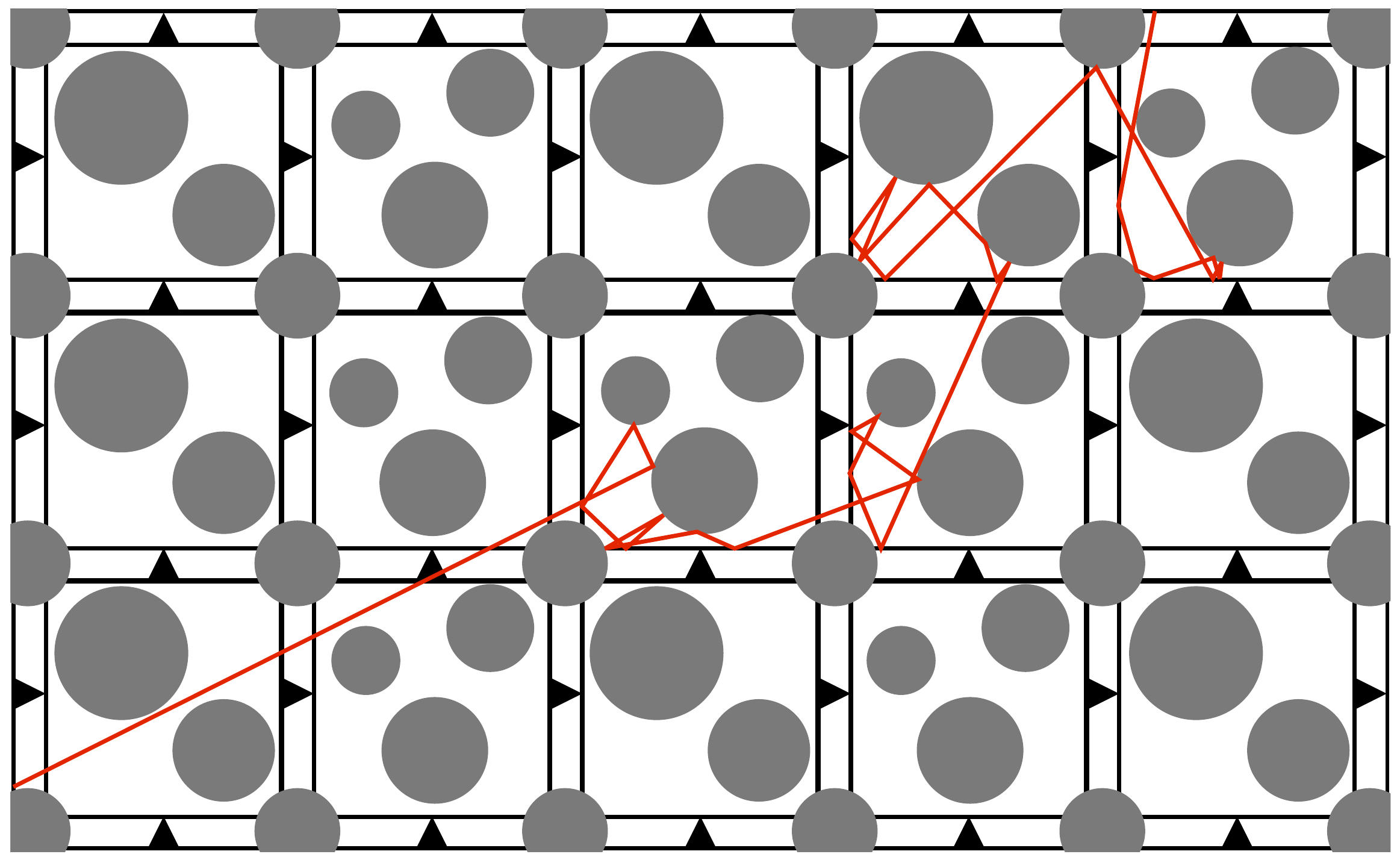} 
   \caption{Schematic diagram illustrating the qualitative features of our two-dimensional model. The medium is composed of two different types of square cells separated by walls which allow particles to pass only in the direction of the arrowheads. The zig-zag line shows a path of a particle through the medium.}
   \label{fig:0}
\end{figure}

Our study concerns an idealization of the aperiodic Lorentz gas with semipermeable walls illustrated in Figure~\ref{fig:0}. In each cell, there is a configuration of scatterers drawn independently from the same probability distribution. In our case, the distribution is Bernoulli, so that there are two possible configurations to choose from inside each cell. As an important aspect, the environment thus obtained is quenched; once the scatterer configurations have been randomly chosen, they are frozen for good, and the only randomness that remains is in the initial data of the particle. Between the cells are semipermeable walls that allow the particle to pass through from left to right and from bottom to top, as shown by the arrowheads, but not in the opposite directions. This model may be thought of as describing the propagation of particles in an exotic, anisotropic, medium. 

Notice that, as a significant simplification, there is no recurrence; once a particle leaves a cell, it never returns to the same cell again. Yet a particle \emph{can} occupy a single cell for an arbitrarily long time before moving on to a neighboring one---albeit a long occupation time has a small probability. Moreover, where, when, and in which direction the particle exits a cell depends heavily on the scatterer configuration inside the cell, in addition to the position and direction of the particle at entry. Inside each cell, the dynamics is chaotic and hyperbolic.

In our one- and two-dimensional idealizations, the billiard dynamics is replaced by discrete dynamical systems acting in each cell. In other words, acting on the particle's current position by a map associated to the current cell gives its position one time unit later. In dimension one the maps associated to the cells are expanding circle maps while in dimension two they are hyperbolic toral automorphisms. Such choice of maps retains the chaotic and hyperbolic nature of the problem. A closely related model has been studied in \cite{AyyerStenlund, AyyerLiveraniStenlund}.

Our objective is to understand certain statistical properties of the motion. More precisely, we are interested in how the particle distribution evolves with time when the initial distribution is uniform and supported on one initial cell. We make several analytical propositions, which we verify numerically. We show that, on the average, the particles follow a linear drift and that, after taking a suitable scaling limit, the cumulative distribution of the fluctuations about the mean is Gaussian. Moreover, the drift and variance are the same for (almost) all environments drawn from the same distribution. Nevertheless, the fine structure of the particle distribution is peculiar due to the quenched environment. In particular, the density function does \emph{not} converge to that of a normal distribution. In fact, it does not converge at all in the aforementioned  scaling limit.

\subsection*{Acknowledgements}
We are indebted to Arvind Ayyer and Joel \linebreak Lebowitz for stimulating discussions. Tapio Simula is supported by the Japan Society for the Promotion of Science Postdoctoral Fellowship for Foreign Researchers. Mikko Stenlund is partially supported by a fellowship from the Academy of Finland.

%%%%%%%%%%%%%%%%%%%%%%%%%
%%%%%%%%%%%%%%%%%%%%%%%%%

\section{One-dimensional model: expanding maps on the circle} 
\subsection{Introduction}
Imagine tiling the nonnegative half line $[0,\infty)$ so that each interval---or tile---$I_k=[k,k+1)$ with $k\in\N$ carries a label $\omega(k)$ that equals either $0$ or $1$. Such a tiling can be realized by flipping a coin for each $k$ and encoding the outcomes in a sequence $\omega=(\omega(0),\omega(1),\dots)\in\{0,1\}^\N$ called the \emph{environment}. The coin could be balanced but the tosses are independent, with $\operatorname{Prob}(\omega(k)=i)=p_i$ for all $k$. Below, $\bP_{p_0}$ will stand for the corresponding probability measure on the space of Bernoulli sequences $\omega$. In each experiment we freeze the environment---meaning that we work with one \emph{fixed} sequence $\omega$ at a time.

The dynamics in our model is generated by the following definitions.
Let $A_0,A_1\in \{2,3,\dots\}$ and define the circle maps $T_i(x)=A_ix \mod 1$. 
An experiment comprises iterating many times the map $\mathbb{R}\times \mathbb{S}^1\to \mathbb{R}\times \mathbb{S}^1:(v,x)\mapsto (v+A_{\omega([v])}x-x,T_{\omega([v])}(x))$, where $[v]$ is the integer part of $v$ and $\omega([v])$ is the corresponding component of $\omega$. The actual initial condition is $(x,x)$, with $x\in[0,1)$, although in the treatment below we will also consider initial conditions of the form $(V+x,x)$ with $V\in\N$ arbitrary. Let $\bfP$ denote the Lebesgue measure (\ie, the uniform probability distribution) on the circle $\mathbb{S}^1$ and $\bfE$ the corresponding expectation.

The difference $A_{\omega([v])}x-x$ above can be thought of as the jump experienced by a particle. In other words, the evolution of the $v$-coordinate \emph{on the line} keeps track of the journey of a particle, from the point of view of the particle. The $x$-coordinate is just the projection onto the circle (and is needed for technical purposes): if $(v_n,x_n)$ is the image of $(v_0,x_0)=(V+x,x)$ after $n$ steps, then $x_n=v_n \mod 1$. 

The model thus describes the deterministic motion of a particle in a randomly chosen, but fixed, environment.  In probability jargon, the particle performs a deterministic walk in a quenched random environment. The challenge is that the \emph{map} determining the position of the particle in the nearest future depends on the tile $I_{[v]}$ the particle is in through the label $\omega([v])$ of the tile. That is, the motion of the particle is guided by the \textit{a priori} chosen environment.

\begin{example}\label{exa:favorite}
A concrete example is obtained by choosing $A_0=2$, $A_1=3$, and $p_0=p_1=\frac{1}{2}$.
\end{example}

%Listing the first few iterates should be useful to illustrate what is going on:
%\begin{align*}
%(v_0,x_0) & = (V+x,\,x) \quad\text{with}\quad 0\leq x<1,\, V\in\N\\
%(v_1,x_1) & = \left(V+x+(A_{\omega(V)}-1)x,\,T_{\omega(V)}(x)\right) \quad\text{since $[x]=0$}\\
%(v_2,x_2) & = \left(V+x+(A_{\omega(V)}-1)x+(A_{\omega([v_1])}-1)T_{\omega(V)}(x) ,\,T_{\omega([v_1])}\circ T_{\omega(V)}(x)\right). 
%\end{align*}
In general, 
\beqn
v_n=[v_0]+x+(A_{\omega([v_0])}-1)x+\sum_{j=1}^{n-1}(A_{\omega([v_{j}])}-1)T_{\omega([v_{j-1}])}\circ\dots\circ T_{\omega([v_0])} (x).
\eeqn
Notice that already in the expression of $v_2$, the initial $V$ and $x$ and the environment $\omega$ are present in a complicated fashion:
\beqn
A_{\omega([v_1])}=A_{\omega([x+(A_{\omega(V)}-1)x])}.
\eeqn
The complication referred to earlier is obvious; for each $x$ (and $\omega$) we have to \emph{compute} which map the symbol $A_{\omega([v_1])}$ stands for. Each $v_n$ ($n\geq 2$) is generically a piecewise affine function of $x$ and the number of discontinuities grows exponentially with $n$.

We anticipate that, for large values of $n$, $v_n$ behaves statistically (in the weak sense) like \
 \beq\label{eq:stat-conjecture}
 v_n \approx \cN(nD,n\sigma^2),
 \eeq
where $D$ is a deterministic number called the \emph{drift} and $\cN(nD,n\sigma^2)$ stands for a real-valued, normally distributed, random variable with mean $nD$ and variance $n\sigma^2$. In principle, $D$ and $\sigma^2$ could depend on the environment $\omega$, but remarkably it turns out that they do not, as long as the environment is \emph{typical}. By typicality we mean that $\omega$ belongs to a set whose $\bP_{p_0}$-probability is one and whose elements enjoy good statistical properties such as the convergence of 
$
\frac1n\#\{k<n\,|\, \omega(k)=0\}
$
to the limit $p_0$.

It is reasonable to expect that the limit
\beqn
\lim_{n\to\infty} \frac{v_n(x)}{n}
\eeqn
exists and has the same value for almost all $x$ \footnote{This cannot hold for all $x$. For instance, if $x=\frac{1}{k A_{\omega(0)}}$, then $1=v_k=v_{k+1}=\dots$ and the process stops.}.
Thus, we are led to conclude that
\beqn
D= \bfE \! \left( \lim_{n\to\infty} \frac{v_n(x)}{n}\right) = \lim_{n\to\infty} \frac{1}{n} \bfE \!\left(v_n(x)\right).
\eeqn
The final equality follows from the bounded convergence theorem.

Let us consider the (asymptotically) centered random variable 
\beqn\label{eq:Xdef}
X_n=v_n-nD,
\eeqn
which measures the fluctuations of $v_n$ relative to the linear drift. Provided \eqref{eq:stat-conjecture} is true, $X_n$ is approximately Gaussian with variance $n\sigma^2$. More precisely, we would like to know if $\frac{1}{\sqrt n}X_n$ converges in distribution to $\cN(0,\sigma^2)$. By definition, this means that, for any fixed $y\in\bR$,
\beqn
\lim_{n\to\infty}\bfP \! \left(\frac{1}{\sqrt n}X_n \leq y\right) = \frac{1}{\sqrt{2\pi}\sigma} \int_{-\infty}^y e^{-s^2/2\sigma^2} \,ds.
\eeqn
%or equivalently, for any fixed $\lambda\in\bR$,
%\beqn
%\lim_{n\to\infty}\bfE \! \left(e^{i\lambda X_n/\sqrt{n}}\right) = e^{-\lambda^2\sigma^2/2}.
%\eeqn 

%%%%%%%%%%%%%%%%%%%%%%%%%%%%%%

\subsection{Markov partition}
We next reduce the deterministic walk in a random environment to a random walk in a random (still quenched) environment which is easier to treat. This can be done using a Markov property of the tiling that allows us, in the statistical sense, to ignore the exact position of the particle and only keep track of  the tile it is occupying.

Let $[\,\cdot\,]$ denote the integer part of a number. If we define
\beq\label{eq:V&x}
V_n=[v_n]\quad\text{and}\quad x_n=v_n-[v_n],
\eeq
then the earlier dynamics with the initial condition $(v_0,x_0)=(V_0+x_0,x_0)$ is equivalent to
\beq\label{eq:recursion}
\begin{split}
V_{n+1} & =V_{n}+\left[A_{\omega(V_{n})}x_{n}\right] \\
x_{n+1} & =A_{\omega(V_{n})}x_{n}-\left[A_{\omega(V_{n})}x_{n}\right].
\end{split}
\eeq

Recall our convention $[0,\infty)=\bigcup_{k=0}^\infty I_k$, where $I_k=\left[k,k+1\right)$ is called a tile. Suppose now that $v_n\in I_k$. This is equivalent to $V_n=k$. As before, we are interested in the probability distribution of $v_n$, but this time only at the level of tiles. More precisely, we wish to know the probability distribution of $V_n$. This is the probability vector  $\rho^{(n)}=(\rho^{(n)}_0,\rho^{(n)}_1,\dots)$ where the numbers
\beqn
\rho^{(n)}_k = \bfP\!\left(V_n=k\right)
\eeqn 
are such that $\sum_{k=0}^\infty \rho^{(n)}_k=1$.
%In particular, 
%\beq\label{eq:init-density}
%\rho^{(0)}=(1,0,0,\dots).
%\eeq

We now consider the dynamical system being initialized with the condition $(v_0,x_0)=(V_0+x_0,x_0)$, where $V_0\in\N$ and $x_0\in[0,1)$ are random variables; $x_0$ is uniformly distributed and independent of $V_0$, but for the moment we do not specify the distribution $\rho^{(0)}_k = \bfP(V_0=k)$ which is the initial distribution of the chain $V_n$. 

Since each $A_i I_k$ is exactly the union of a few of the intervals $I_{k'}$\footnote{$A_i$ maps $\left[k+\frac{l}{A_i},k+\frac{l+1}{A_i}\right)$ affinely onto $[k+l,k+l+1)$.}, the collection $\{I_k\}$ is a simultaneous Markov partition for the two maps. We then obtain the Markov property
\beqn
\bfP(V_{n}=k_n\,|\, V_{n-1}=k_{n-1},\dots,V_0=k_0) \equiv \bfP(V_{n}=k_n\,|\,V_{n-1}=k_{n-1})
\eeqn 
for arbitrary histories. Thus, $V_n$ is a time-homogeneous Markov chain on the countably infinite state space $\N$ with the transition probabilities
\beqn
\gamma_{k\rightarrow k+l} =
\begin{cases}
\bfP(v_{n+1}\in I_{k+l}\,|\,v_n\in I_k) = \frac{1}{A_{\omega(k)}} & \text{if}\quad l\in\{0,1,\dots,A_{\omega(k)}-1\}, \\
0 & \text{otherwise}.
\end{cases}
\eeqn
Notice that the above holds for any environment, $\omega$, but the resulting Markov chain does depend on the choice of $\omega$.

Defining the transition matrix
$
\Gamma = {(\gamma_{k\to k'})}_{k,k'},
$
$\rho^{(n)}=\rho^{(n-1)}\Gamma$. Thus,
\beq\label{eq:density}
\rho^{(n)} = \rho^{(0)} \Gamma^n
\eeq
for an arbitrary initial distribution. 
In particular,
\beqn
\rho^{(0)}=(1,0,0,\dots)
\eeqn
corresponds to initializing the dynamical system at $(v_0,x_0)=(x,x)$ with $0\leq x<1$ arbitrary.

In principle, \eqref{eq:density} provides us with complete statistical understanding of the dynamics. For instance, the drift can be expressed as
\beqn
D=\lim_{n\to\infty} \frac{1}{n}\,\bfE(v_n)=\lim_{n\to\infty} \frac{1}{n}\,\bfE(V_n) = \lim_{n\to\infty} \frac{1}{n}\sum_{k=0}^\infty k\rho^{(n)}_k.
\eeqn
In practice, calculating $\Gamma^n$ for large values of $n$ is difficult.

%%%%%%%%%%%%%%%%%%%%%%%%%

\subsection{Drift and variance}
For each $(i,j)\in\{0,1\}^2$ the transition probability at time $n$ from a tile labeled $i$ to a tile labeled $j$ is 
\beqn
\alpha_{ij}(n)=\bfP(\omega(V_{n+1})=j\,|\,\omega(V_n)=i).
\eeqn
The analysis of this quantity is subtle, because it depends on the tiling. For instance, if $\omega=(0,0,\dots)$, then $\bfP(\omega(V_{n+1})=0\,|\,\omega(V_n)=0)=1$. 

The conditional probability $\bfP\times \bP_{p_0}(\omega(V_{1})=j\,|\,\omega(V_0)=i)$ equals
\beqn
\alpha_{ij}^* =\delta_{ij}\left(\frac{1}{A_i}+\left(1-\frac{1}{A_i}\right)p_i\right)+(1-\delta_{ij})\left(1-\frac{1}{A_i}\right)p_j,
\eeqn
because the elements $\omega(k)$ of the tiling are independent. Here $p_i$ is the Bernoulli probability of getting an $i$ in the tiling. We think of $\alpha_{ij}^*$ as an effective transition probability which only depends on the statistical properties of the tiling.

We expect the actual transition probability $\alpha_{ij}(n)$ to converge to the effective value $\alpha_{ij}^*$ with increasing time,
\beqn
\lim_{n\to\infty}{\alpha_{ij}(n)}=\alpha_{ij}^*.
\eeqn
This is so because, as $n$ increases, the position of the particle at time $n$ depends on the tiling on an increasing subinterval of $[0,\infty)$ and should therefore reflect increasingly the statistics of the tiling instead of its local details. 

Moreover,
\beq\label{eq:p-matrix}
\lim_{n\to\infty}(\alpha^*)^n = 
\begin{pmatrix}
p & 1-p \\
p & 1-p
\end{pmatrix}
\eeq
for a $p\in(0,1)$ that can be found by diagonalizing $\alpha^*$ or by solving the equilibrium equation $(p,1-p)\alpha^*=(p,1-p)$:
\beqn
p=\frac{p_0\left(1-\frac{1}{A_1}\right)}{1-p_1\frac{1}{A_0}-p_0\frac{1}{A_1}}=\frac{p_0A_0(A_1-1)}{A_0A_1-p_1A_1-p_0A_0}.
\eeqn
For instance, in the case of Example~\ref{exa:favorite} we obtain
$
p=
\frac{4}{7}. 
$

Notice that, for any probability vector $(q,1-q)$,
\beqn
(q,1-q)\lim_{n\to\infty}(\alpha^*)^n=(p,1-p).
\eeqn
The probability vector
\beqn
(q,1-q)\prod_{n=0}^k \alpha(n)=(q(k),1-q(k))
\eeqn
will converge to some $(q^*,1-q^*)$, because $\alpha(n)\to\alpha^*$.
In fact,
\beqn
\begin{split}
\lim_{N\to\infty} (q,1-q)\prod_{n=0}^{2N} \alpha(n) &=\lim_{N\to\infty}  (q(N),1-q(N))\prod_{n=N+1}^{2N} \alpha(n)
\\
&= (q^*,1-q^*) \lim_{N\to\infty} (\alpha^*)^N = (p,1-p),
\end{split}
\eeqn
as $N\to\infty$.

We interpret the result above so that $\bfP(\omega(V_n)=0)\to p$ and $\bfP(\omega(V_n)=1)\to 1-p$ as $n\to\infty$. That is, along a given (typical) trajectory, the fraction of time the particle spends in a tile labeled 0 is $p$:
\beq\label{eq:p}
\lim_{n\to\infty}\frac{\#\{k<n\,|\, \omega(V_k)=0\}}{n}=p.
\eeq
Notice that $p$ does not depend on the (typical) tiling.

%%%%%%%%%%%%%%%%%%%%%%%%%%%

\subsubsection{Drift}\label{subsubsec:drift}
Define the jumps $\xi_i=V_i-V_{i-1}$ ($i\geq 1$). Then $V_n=\sum_{i=1}^n \xi_i$. We also denote $\xi^{(j)}$ a random variable that takes values in $\{0,\dots,A_j-1\}$ with uniform distribution. For a (typical) tiling,
\beqn
\begin{split} 
D  = \lim_{n\to\infty}\frac{\bfE(V_n)}{n} & = \lim_{n\to\infty}\frac{\sum_{i=1}^n\bfE(\xi_i)}{n} =  p\bfE(\xi^{(0)})+(1-p)\bfE(\xi^{(1)})
\\
& = p\frac{A_0-1}{2}+(1-p)\frac{A_1-1}{2} = \frac{pA_0+(1-p)A_1-1}{2}.
\end{split}
\eeqn
We also claim that this equals $\lim_{n\to\infty}\frac{V_n}{n}$ for almost all $x$. In the case of Example~\ref{exa:favorite}, $D=\frac{5}{7}$.

%%%%%%%%%%%%%%%%%%%%%%%%%%%%

\subsubsection{Variance}\label{subsubsec:variance}
The variance is 
\beqn
\sigma^2=\lim_{n\to\infty}\Var\!\left(\frac{X_n}{\sqrt n}\right)=\lim_{n\to\infty}\frac{\Var(V_n)}{n}.
\eeqn

Let us assume $A_0\leq A_1$ and study the process $W_n=\sum_{i=1}^n\zeta_i$ having the i.i.d.\@ increments $\zeta_i$ whose distribution is $\operatorname{Prob}(\zeta_1=k)=\frac{p}{A_0}+\frac{1-p}{A_1}$ if $0\leq k<A_0$ and $\operatorname{Prob}(\zeta_1=k)=\frac{1-p}{A_1}$ if $A_0\leq k<A_1$.
The increments have been chosen so that $W_n$ mimics $V_n$ as closely as possible and has the same variance. For instance, staying in the same tile ($\zeta_1=0$) has probability $\frac{p}{A_0}+\frac{1-p}{A_1}$, where $p$ is the probability of being in a tile labeled $0$ and $\frac{1}{A_0}$ is the probability of staying in that tile, while the second term accounts similarly for the case of label 1.
Then $\operatorname{Mean}(W_n)=n\operatorname{Mean}(\xi_1)=nD$. Setting $K(m)=\sum_{k=0}^{m-1} k^2 = \frac 13 m^3 + \frac 12 m^2 + \frac 16 m$, the variance of $W_n$ is
\beqn
\begin{split}
\Var(W_n) & = n\Var(\zeta_1)= n\left(\frac{p}{A_0}K(A_0)+\frac{1-p}{A_1}K(A_1)-D\right).
\end{split}
\eeqn

Using the values $p=\frac{4}{7}$ and $D=\frac{5}{7}$ obtained from Example~\ref{exa:favorite}, the formula above gives $\frac{1}{n}\Var(W_n)=\frac{24}{49}$.

%%%%%%%%%%%%%%%%%%%%%

\subsection{Sensitivity on the initial condition} 
Let us next consider the Lyapunov exponent
\beqn
\lambda=\lim_{n\to\infty}\frac{1}{n}\ln\frac{dv_n}{dx} = \lim_{n\to\infty}\frac{1}{n}\sum_{k=1}^n\ln\frac{dv_k}{dv_{k-1}}
\eeqn
which measures the exponential rate at which two nearby initial points drift apart under the dynamics. Above, the chain rule has been used. Recall the notation introduced in \eqref{eq:V&x} and that $v_0=x_0=x$. As $v_{k}=v_{k-1}+\left(A_{\omega(V_{k-1})}-1\right)x_{k-1}$,
\beqn
\frac{dv_k}{dv_{k-1}}=1+\frac{d A_{\omega(V_{k-1})}}{dv_{k-1}}x_{k-1}+\left(A_{\omega(V_{k-1})}-1\right)\frac{dx_{k-1}}{dv_{k-1}}.
\eeqn
With probability zero $v_{k-1}$ is an integer, in which case $v_{k-1}=V_{k-1}$, $x_{k-1}=0$, and the process stops. We assume that $v_{k-1}$ is not an integer.
Then $\frac{d A_{\omega(V_{k-1})}}{dv_{k-1}}=0$, $\frac{dx_{k-1}}{dv_{k-1}}=1$, and
$
\frac{dv_k}{dv_{k-1}}=A_{\omega(V_{k-1})}
$,
such that, by \eqref{eq:p}, 
\beqn
\lambda=\lim_{n\to\infty}\frac{1}{n}\sum_{k=1}^n\ln A_{\omega(V_{k-1})} = p\ln A_0+(1-p)\ln A_1 
\eeqn
and is positive. Roughly speaking, the distance between two very nearby trajectories thus grows like $e^{\lambda n}=\left(A_0^pA_1^{1-p}\right)^n$, which is tantamount to chaos.

%%%%%%%%%%%%%%%%%%%%%%

\subsection{Numerical study}
% TILING
In order to study the model introduced above numerically we first create the random tiling (or environment) $\omega$ of length $2n+1$ where $n$ is the number of jumps to be performed in a single trajectory. This guarantees that every possible path fits inside the tiling although some computational effort could be saved by choosing the number of tiles closer to $[nD]$. Notice also that while in principle new tiles could be added dynamically to the end of the tiling as required, it is computationally far more efficient to construct the tiling as a static entity in the beginning of the computation.  

In practice, the tiling is generated by producing a vector of pseudo-random numbers distributed uniformly on the interval $(0, 1)$ using the Mersenne Twister algorithm. The label of the tile $\omega(k)$ is then obtained by rounding the number on each tile $k$ to the nearest integer. The computational tiling $\hat{\omega}$ is finalized by the operation $\hat{\omega}(k) = \omega(k) (A_1 - A_0) + A_0$ yielding a vector whose each element is either 2 or 3. %The first 21 elements in our specific case are $\hat{\omega}(k=[1,21])=\{332332233323323223333\}$. 

% TRAJECTORIES
Each ensemble member (particle trajectory) is initialized by generating a pseudo-random number to determine the starting point $x_0\in(0,1)$ of the trajectory. The subsequent particle positions are determined by the underlying tiling. The jumping process could be performed deterministically by keeping track of the exact position $v_n$ of the particle. However, the Markov property of the process provides us a superior way of obtaining the desired statistics stochastically. In this algorithm, before every jump, we sample a new pseudo-random number $d$ from the interval $(0,\hat{\omega}(k))$ depending on the current tile $k$. Then a jump to the tile $k+[d]$ is made and the whole procedure is repeated $n$ times to produce a single trajectory. 

\begin{figure}
   \centering
   \includegraphics[width=0.7\linewidth]{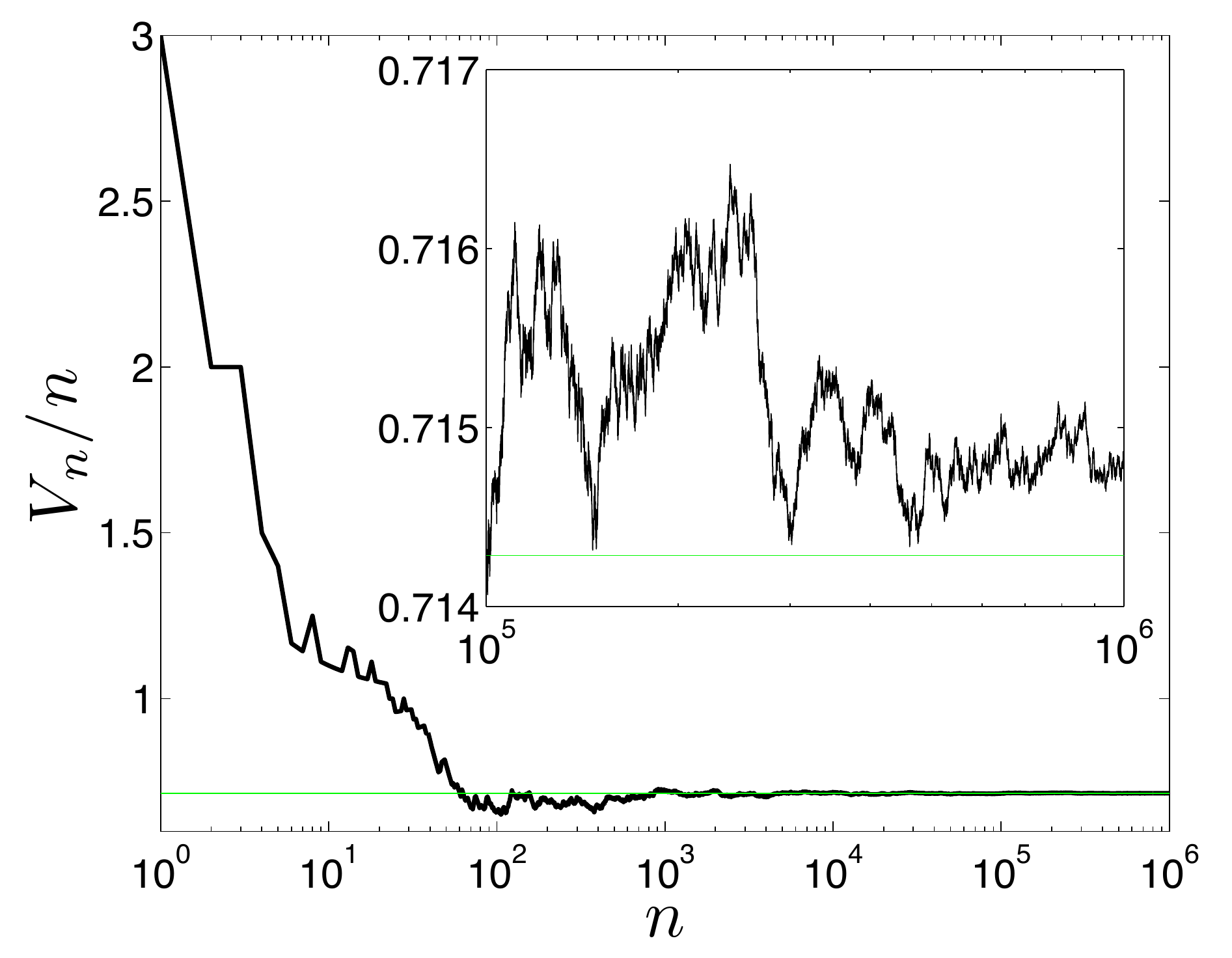} 
	   \caption{Position of a particle $V_n$ divided by the number of jumps $n$ taken for a single typical trajectory. The horizontal line denotes the exact value for the drift $D=5/7$. The inset shows a blow-up of a part of the main figure.}
   \label{fig:1}
\end{figure}

% RESULTS
Figure~\ref{fig:1} shows a typical trajectory of $n=10^6$ jumps obtained using the above prescription. The position of the particle $V_n$ divided by the number of jumps $n$ taken is clearly seen to saturate to the analytical value for the drift $D=5/7$, plotted as a straight line in the figure. The inset shows the late-time evolution of the drift of the particle.

\begin{figure}
   \centering
   \includegraphics[width=0.7\linewidth]{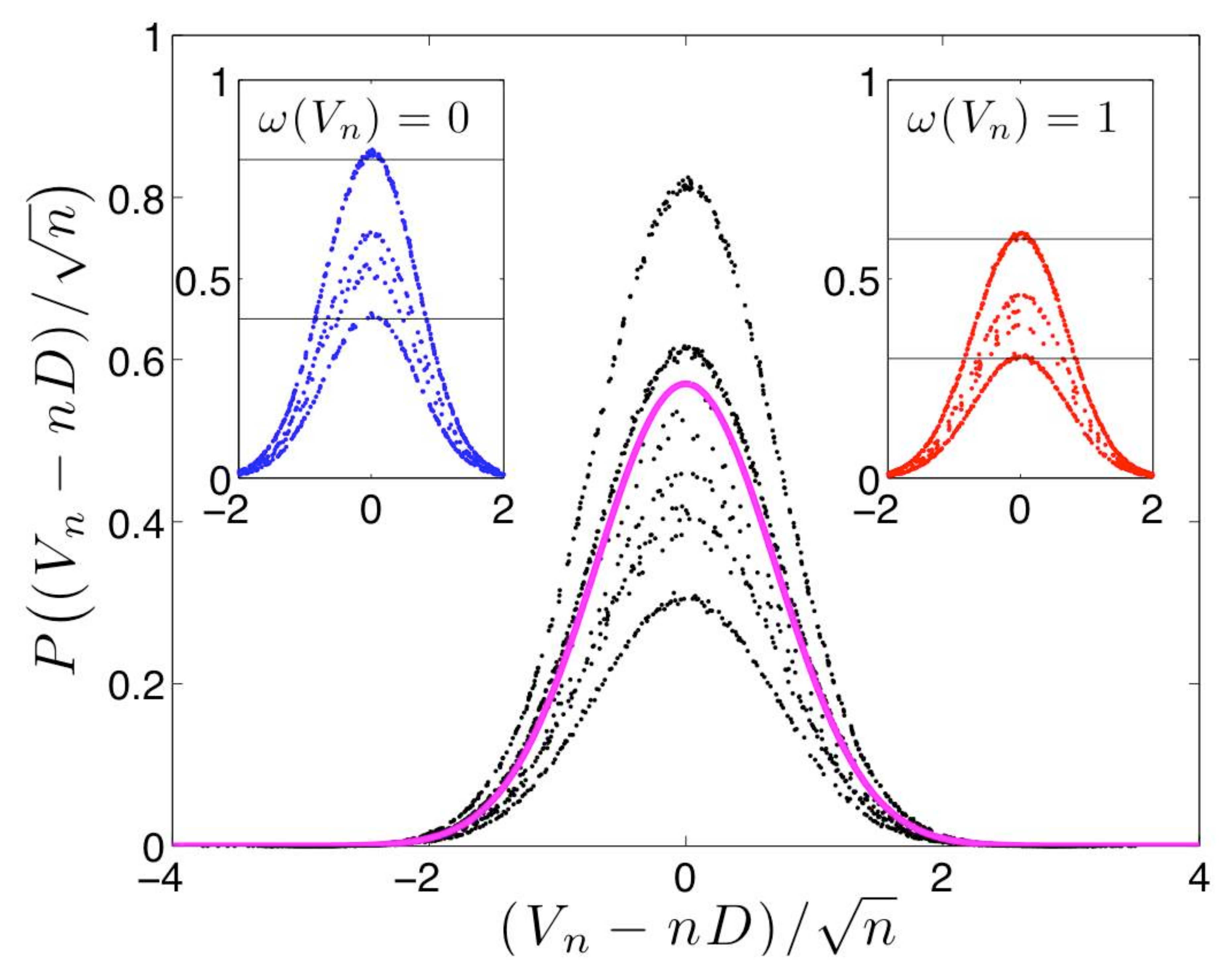}%{Kuva3Fixar.pdf} 
   \caption{Probability density of the random variable $(V_n-nD)/\sqrt{n}$. The solid curve is the  Gaussian with zero mean and variance $\sigma^2=24/49$. The insets in the left- and right-hand sides show the probability densities for the subsets of trajectories ending to a tile labeled by 0 and 1, respectively. The horizontal lines at levels 0.4, 0.8, 0.3 and 0.6 in the insets are plotted to guide the eye.}
   \label{fig:3}
\end{figure}

Figure~\ref{fig:3} displays the computed probability density function for the random variable $(V_n-nD)/\sqrt{n}$ obtained using $10^7$ trajectories of length $n=10^5$. The solid curve is a normalized Gaussian function with zero mean and variance of $\sigma^2=24/49$. Each point in the figure indicates the relative frequency that a trajectory ends in the corresponding tile after $n$ steps. In the left- and right-hand side insets only trajectories ending in tiles labeled by 0 and 1, respectively, are considered. If the graph in the right-hand side inset is vertically stretched by the factor $\frac{p}{1-p}$, it becomes practically overlapping with the one in the left-hand side inset. In the figure one can discern several Gaussian shapes, all of which are very well approximated by the analytical Gaussian after normalization with a suitable constant.

\begin{figure}
   \centering
   \includegraphics[width=0.7\linewidth]{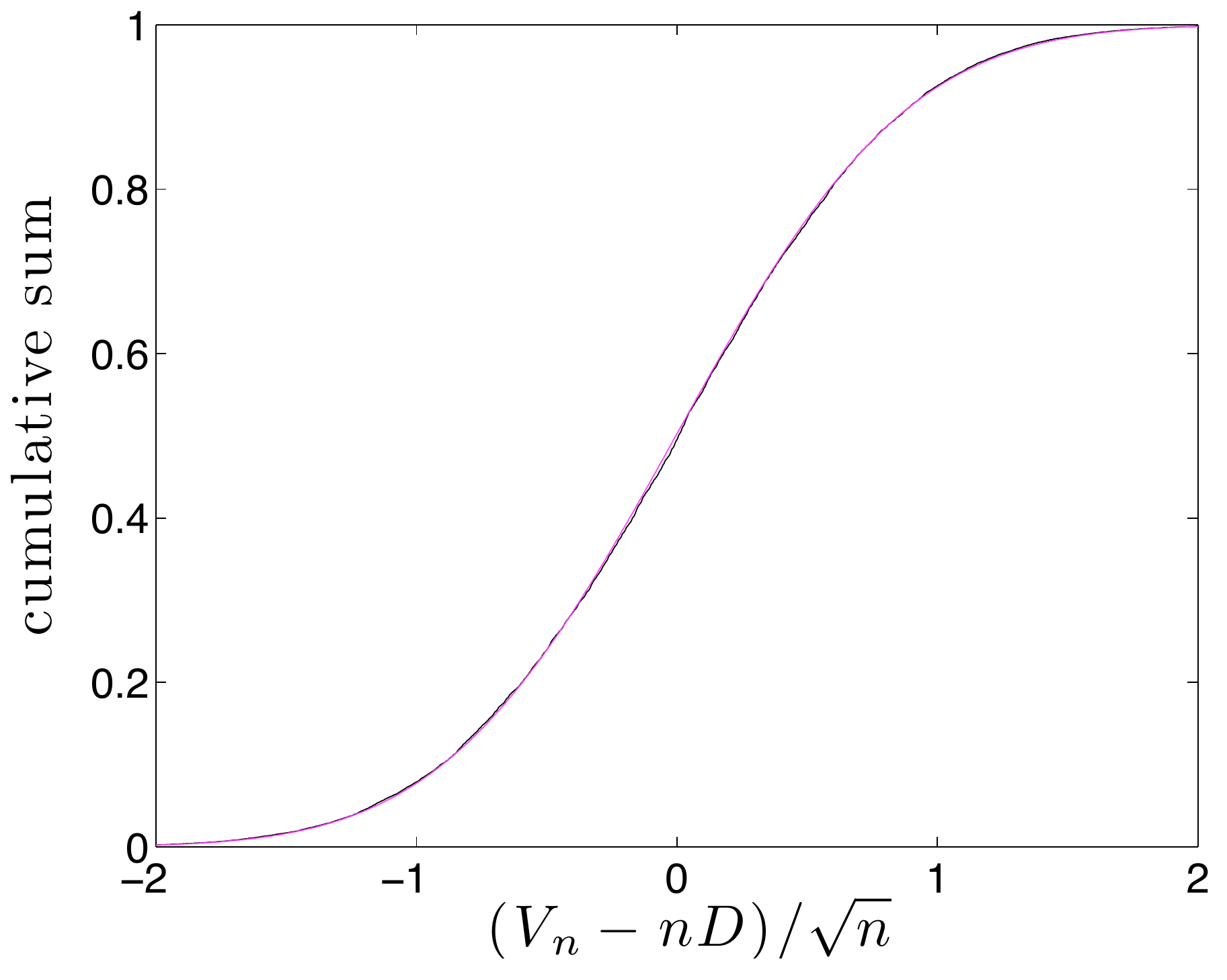} 
   \caption{Cumulative distribution function corresponding to the data shown in Figure~\ref{fig:3}. The two curves shown (computational and analytical) are overlapping.}
   \label{fig:2}
\end{figure}

Figure~\ref{fig:2} shows two cumulative distribution functions, obtained by integrating the numerical and analytical probability densities shown in Figure~\ref{fig:3}. They match to a great accuracy and we are lead to believe that the random variable $X_n/\sqrt{n}=(v_n-nD)/\sqrt{n}$ is, indeed, normally distributed with zero mean and variance $\sigma^2$. We have also analyzed the characteristic function which leads to the same conclusion.

To conclude, our numerical data strongly supports the theoretical analysis presented earlier. Within the numerical accuracy, the distribution is Gaussian with the drift and variance predicted by our analytical calculations. We have performed these numerical experiments using different (fixed) tilings and filling probabilities, $p_i$, and have always arrived at the same conclusion.

%%%%%%%%%%%%%%%%%%%%%%%%%%%%
%%%%%%%%%%%%%%%%%%%%%%%%%%%%

\section{Two-dimensional model: hyperbolic toral automorphisms}
\subsection{Introduction}
We begin by tiling the first quadrant of the plane by unit squares, attaching the label 0 or 1 to each tile. That is, corresponding to each vector $k=(k_1,k_2)\in \N^2$ the tile $[k_1,k_1+1)\times[k_2,k_2+1)$ carries a label $\omega(k)\in\{0,1\}$. The fixed tiling $\omega=(\omega(k))_{k\in \N^2}$ is our environment and $\bP_{p_0}$ stands for the Bernoulli probability measure on the space of such tilings.

The process $v_n$ takes place on the plane and each $A_{i}$ ($i=0,1$) is a matrix with positive integer entries and determinant 1. Such a matrix is hyperbolic, with two eigenvalues, $\lambda>1$ and $\lambda^{-1}$, and the eigenvector corresponding to $\lambda$ points into the first quadrant.
 The formula $T_ix=A_ix\mod 1$ defines a hyperbolic toral automorphism. A precise description of the dynamics is given by the map $\mathbb{R}^2\times \mathbb{T}^2\to \mathbb{R}^2\times \mathbb{T}^2:(v,x)\mapsto (v+A_{\omega([v])}x-x,T_{\omega([v])}(x))$, where $[v]$ is the integer part of $v$. The initial condition is $(x,x)$, with $x\in[0,1)^2$. Let $\bfP$ denote the Lebesgue measure (\ie, the uniform probability distribution) on the torus $\mathbb{T}^2$ and $\bfE$ the corresponding expectation.
 
\begin{example}\label{exa:favorite2}
A concrete example is obtained by choosing $A_0=\left(\begin{smallmatrix} 2 & 1 \\ 1 & 1 \end{smallmatrix}\right)$, $A_1=\left(\begin{smallmatrix} 3 & 1 \\ 2 & 1 \end{smallmatrix}\right)$, and $p_0=p_1=\frac{1}{2}$.
\end{example} 
 
We claim that the limit
 \beqn
D= \bfE \! \left( \lim_{n\to\infty} \frac{v_n(x)}{n}\right) = \lim_{n\to\infty} \frac{1}{n} \bfE\!\left( v_n(x)\right),
\eeqn
called the drift, exists and that the (asymptotically) centered random vector 
\beqn\label{eq:Zdef}
Z_n=(X_n,Y_n)=v_n-nD,
\eeqn
which measures the fluctuations of $v_n$ relative to the linear drift, is approximately Gaussian with covariance matrix $n\sigma^2$. More precisely, $\frac{1}{\sqrt n}Z_n$ converges in distribution to $\cN(0,\sigma^2)$, where $\sigma^2$ is given by \linebreak $\lim_{n\to\infty}\operatorname{Cov}\!\left(\frac{1}{\sqrt n}Z_n,\frac{1}{\sqrt n}Z_n\right)$: denoting $E_{z}=(-\infty,z_1]\times(-\infty,z_2]$ for any fixed $z=(z_1,z_2)\in\bR^2$,
\beqn
\lim_{n\to\infty}\bfP \! \left(\frac{1}{\sqrt n}X_n \leq z_1,\frac{1}{\sqrt n}Y_n\leq z_2 \right) = \frac{1}{2\pi \sqrt{\det \sigma^2}} \int_{E_{z}} e^{-\frac12 s\cdot (\sigma^2)^{-1} s} \,d^2s.
\eeqn

In contrast with the one-dimensional case, the tiling is \emph{not} a Markov partition for the maps, which considerably complicates the analysis of the model.

%%%%%%%%%%%%%%%%%%%%%%%%%%%

\subsection{Drift}
Let us continue to denote $V_n=[v_n]$. We conjecture that the drift vector $D=\left(\begin{smallmatrix}d_1\\ d_2\end{smallmatrix}\right)$ is given by
\beqn
D=pD_0+(1-p)D_1,
\eeqn
where $D_i=\bfE(A_ix-x)=(A_i-\one)\left(\begin{smallmatrix}\frac12\\\frac12\end{smallmatrix}\right)$ equals the average jump under the action of the matrix $A_i$ and $p$ is as in \eqref{eq:p}. The value of $p$ is obtained, as above \eqref{eq:p-matrix}, from an effective transition matrix $\alpha^*$. Its general element $\alpha^*_{ij}$ is the conditional probability $\bfP\times\bP_{p_0}(\omega(V_{1})=j\,|\,\omega(V_0)=i)=\bfP\times\bP_{p_0}(\omega([A_ix])=j\,|\,\omega(0,0)=i)$---the probability of jumping to a tile labeled $j$ when the initial tile is labeled $i$ and when the choice of the tiling is being averaged out.

In practice, $\alpha^*_{ij}$ is computed as follows. We assume that the initial tile is labeled $i$, \textit{i.e.}, $\omega(0,0)=i$. The image of the unit square under $A_i$ is a parallelogram of area one that overlaps with various tiles. The area of intersection of the parallelogram with a tile represents the probability of jumping to that tile. $\alpha^*_{ij}$ can then be computed recalling that each tile is labeled $0$ with probability $p_0$ independently of the others. In the case of Example~\ref{exa:favorite2}, we obtain $D_{0}=\left(\begin{smallmatrix}1\\\frac12\end{smallmatrix}\right)$, $D_1=\left(\begin{smallmatrix}\frac32\\ 1\end{smallmatrix}\right)$, and $\alpha^*=\left(\begin{smallmatrix} \frac14+\frac34\frac12 & \frac34\frac 12 \\ \frac 56\frac12 & \frac16 + \frac 56\frac 12 \end{smallmatrix}\right)=\left(\begin{smallmatrix} \frac58 & \frac38 \\ \frac 5{12} & \frac7{12} \end{smallmatrix}\right)$, which results in $p=\frac{10}{19}$ and $D=\left(\begin{smallmatrix}\frac{47}{38}\\ \frac{14}{19}\end{smallmatrix}\right)$.

%%%%%%%%%%%%%%%%%%%%%%%%%%

\subsection{Numerical study}
% MODEL
As mentioned earlier, the Markov property deployed in the numerical study of the one-dimensional problem where we used a stochastic jumping algorithm does not, unfortunately, apply in the two-dimensional case. Instead, we are forced to compute the particle trajectories fully deterministically which renders the numerical problem difficult. Due to the chaotic nature of the process, the position $v_n$ of the particle must now be represented with an accuracy to approximately $2n$ decimal places in order to keep the accumulation of the numerical rounding errors bounded. This must be done using a software implementation since the double precision float native to the hardware only contains 15 decimal places. 

% TILING AND TRAJECTORIES
We first create the tiling $\omega$ as in the one-dimensional case with the exception that it is now a two-dimensional object. We then choose randomly the initial position of the particle within the unit square. The label $\omega(0,0)$ of the initial tile is then read and the new position of the particle is computed by applying the corresponding map $T_{\omega([v])}$. This jumping procedure is repeated $n$ times. The time to compute a single trajectory increases dramatically as the path length $n$ is increased due to the corresponding increase in the required accuracy of the representation of the position of the particle. 
	
\begin{figure}
   \centering
   \includegraphics[width=0.7\linewidth]{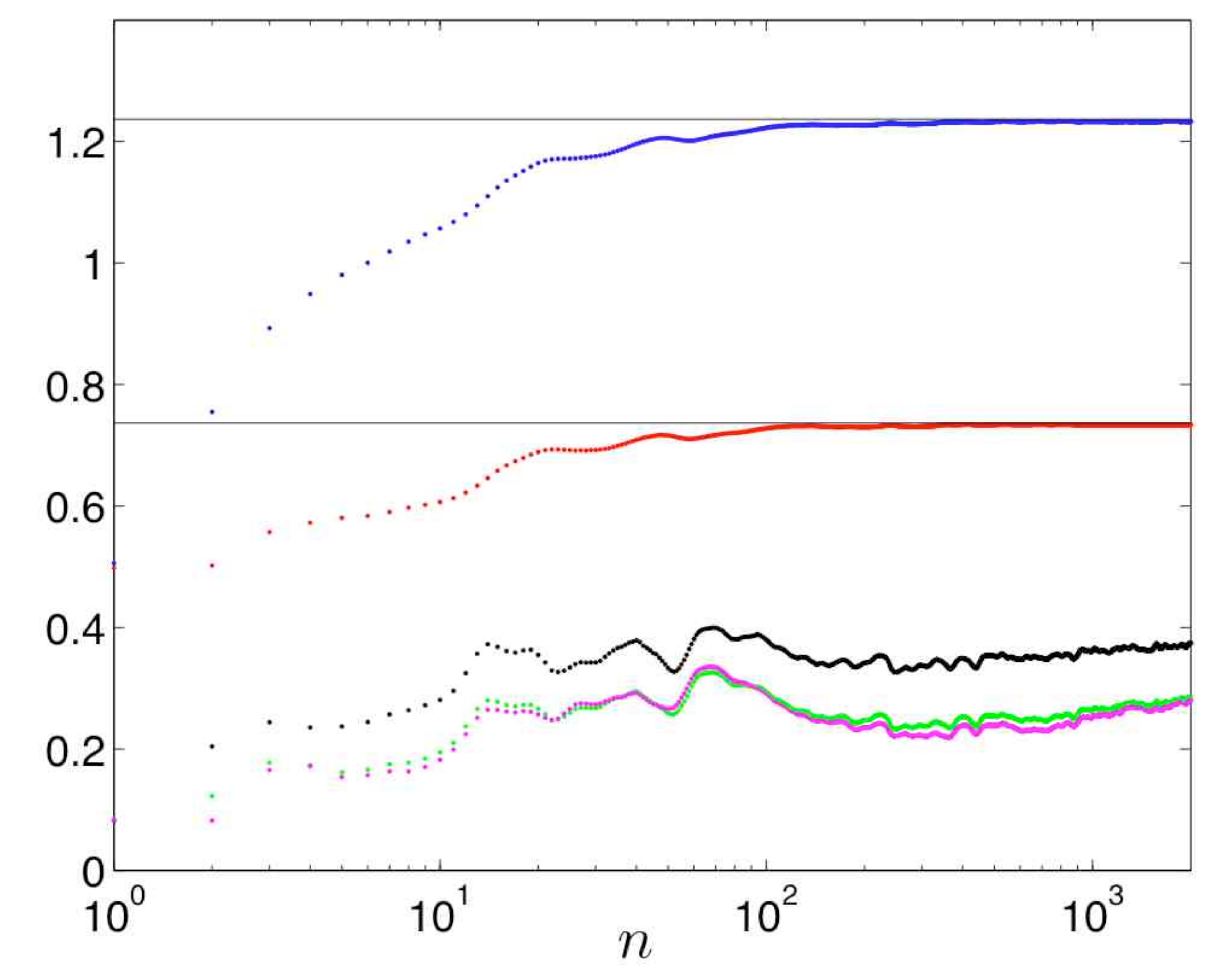} 
   \caption{Values for the $x$ and $y$ components of the drift and the covariance matrix elements $\sigma^2_{11}$, $\sigma^2_{12}=\sigma^2_{21}$, and $\sigma^2_{22}$ as a function of $n$, respectively, in decending order at $n=1000$. The straight lines indicate the analytical values of drift components $d_1$ and $d_1$.}
   \label{fig:4}
\end{figure}

% RESULTS
Figure~\ref{fig:4} shows the convergence of the drifts $D_i$ and that of the covariance matrix elements $\sigma^2_{ij}$. The values in decending order at $n=10^3$ are $d_1$, $d_2$ , $\sigma^2_{11}$, $\sigma^2_{12}=\sigma^2_{21}$, and $\sigma^2_{22}$. The straight lines indicate the analytical values for the drift components. Each data point is composed using $10^4$ trajectories.

\begin{figure}
   \centering
   \includegraphics[width=0.7\linewidth]{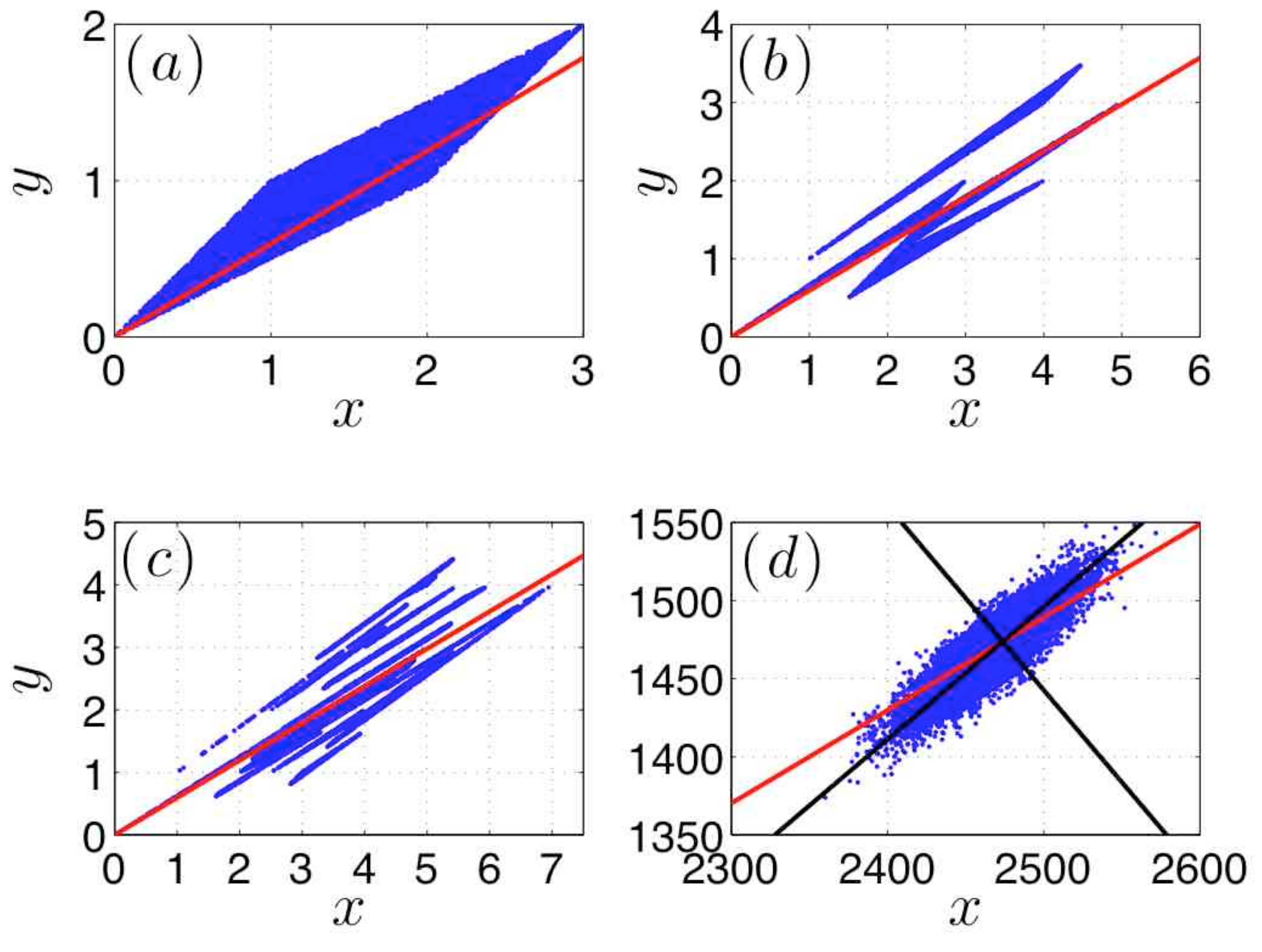} 
   \caption{End points $v_n$ of $10^4$ trajectories in the plane after $n=1$ (a), $n=2$ (b), $n=3$ (c), and $n=2000$ (d) jumps. The straight diagonal lines trace the the drift vector and the cross in frame (d) shows the eigendirections of the covariance matrix. Each frame comprises 10000 data points.}
   \label{fig:7}
\end{figure}

\begin{figure}
   \centering
   \includegraphics[width=0.7\linewidth]{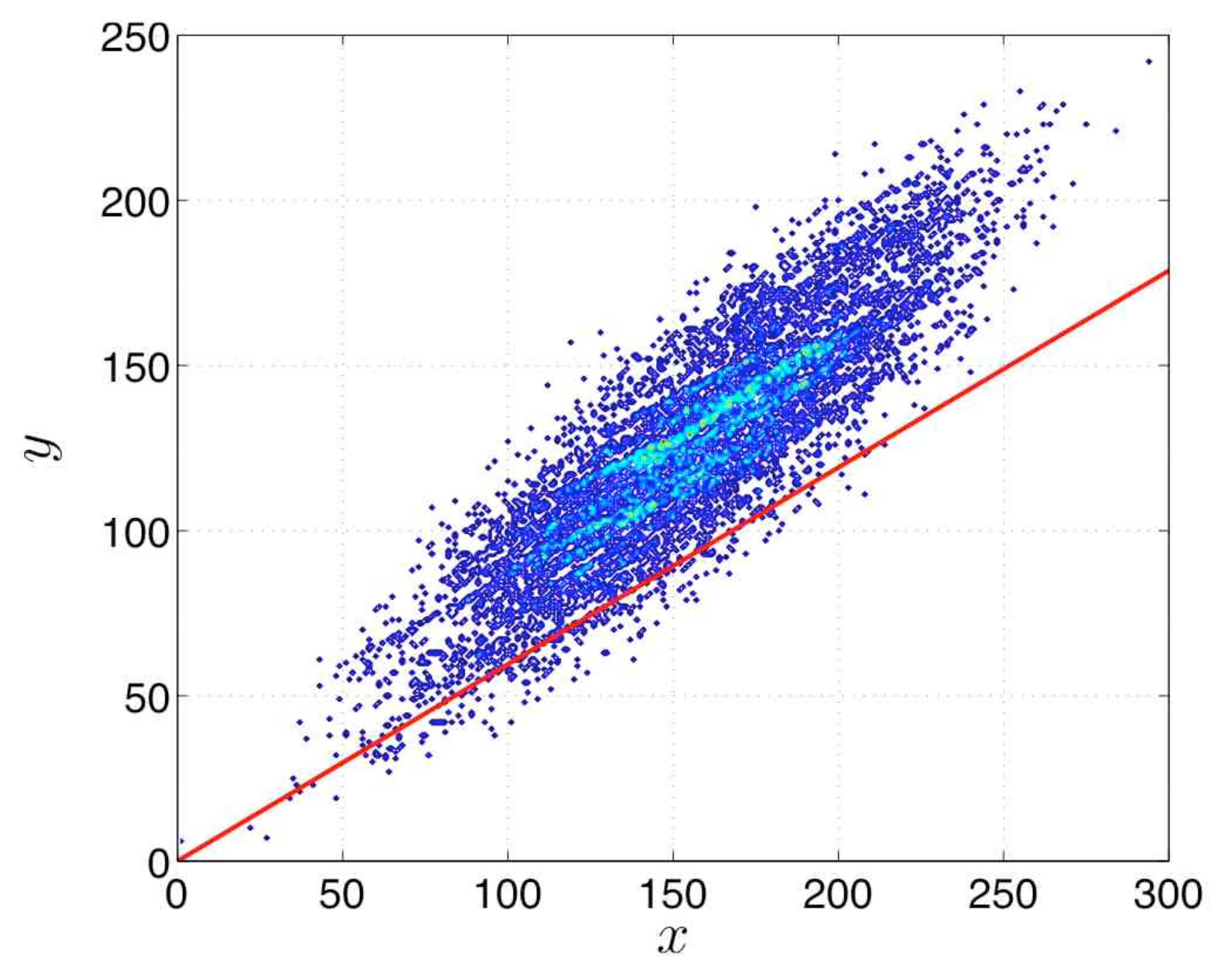} 
   \caption{Contour plot of the particle distribution in the plane after $n=100$ jumps. The straight line shows the direction of the drift vector. }
   \label{fig:8}
\end{figure}

In Figure~\ref{fig:7} (a)--(c) we have plotted the particle positions in the plane after $n=1,2,3$, and $2000$ jumps, respectively. The trajectories were initiated randomly from the unit square. The straight diagonal lines indicate the direction of the drift and the cross in frame (d) denotes the directions of the eigenvectors of the covariance matrix. Since the initial tile in our environment had the label $\omega(0,0)=1$, the frame (a) simply shows how $A_1$ maps the unit square. The subsequent jumps shred the distribution, as illustrated by the frames (b) and (c), because particles in different tiles undergo different transformations. Figure~\ref{fig:8} shows a contour plot of the particle distribution after $n=100$ jumps and reveals prominent  stripes, due to the shredding, which are roughly aligned with the direction of the drift vector.

\begin{figure}
   \centering
   \includegraphics[width=0.7\linewidth]{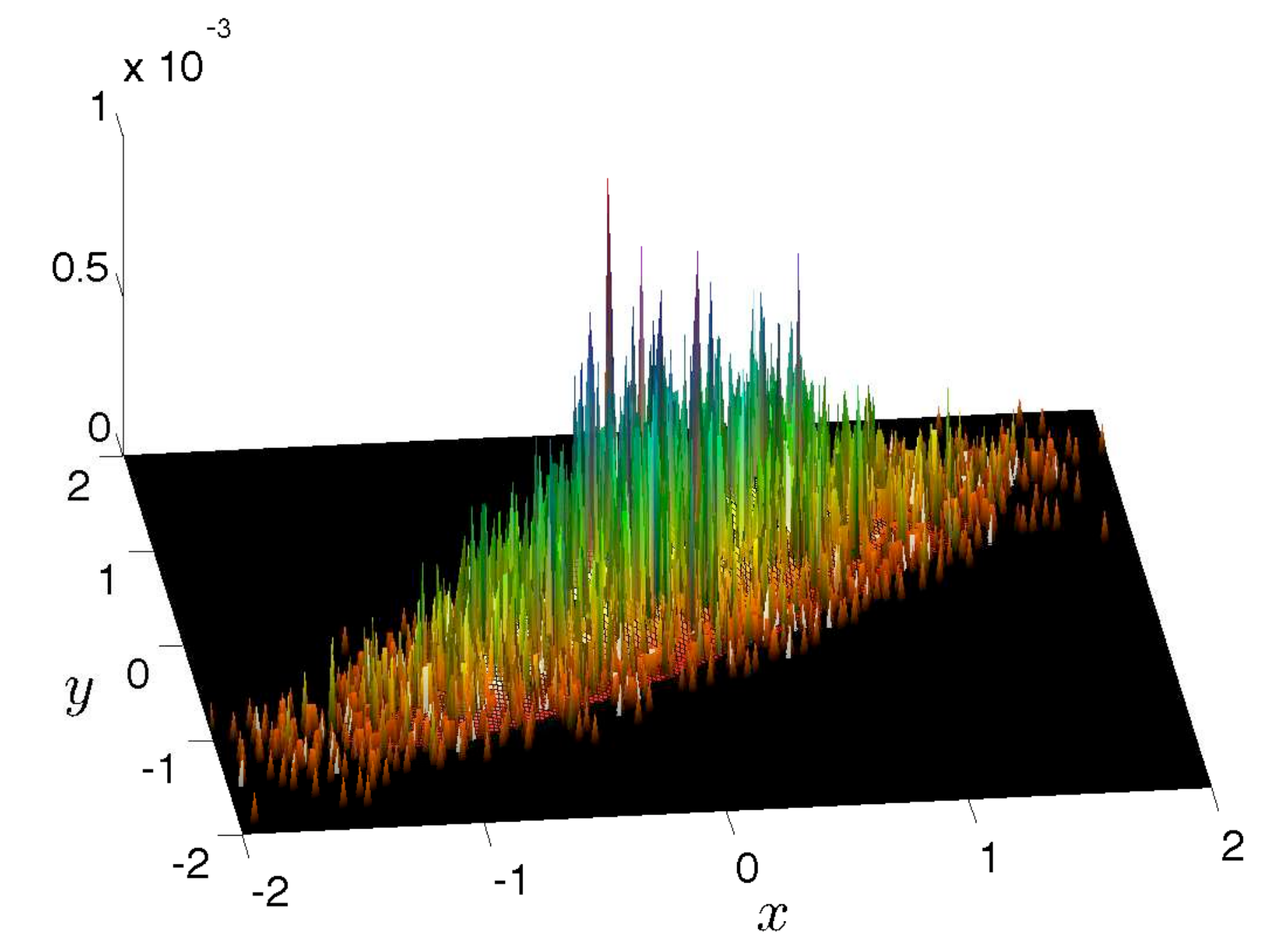} 
   \caption{Probability density function for the random variable $Z_n/\sqrt{n}$. The spikes are an inherent feature of the distribution and the density does not converge to any function. Embedded is the analytical Gaussian function.}
   \label{fig:6}
\end{figure}

\begin{figure}
   \centering
   \includegraphics[width=0.7\linewidth]{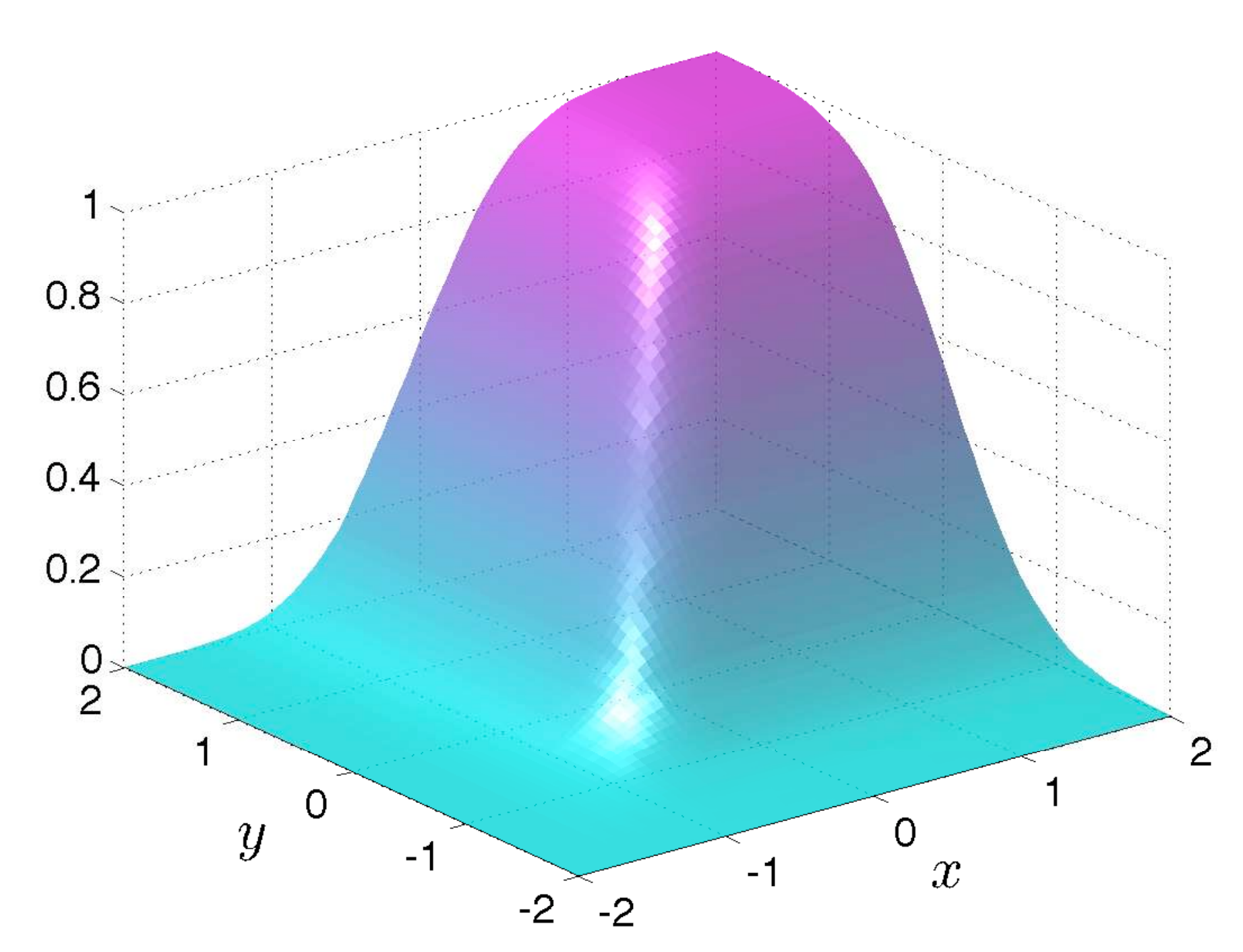} 
   \caption{Cumulative distribution function corresponding to the data shown in Figure~\ref{fig:6}.}
   \label{fig:5}
\end{figure}

Figure~\ref{fig:6} shows the probability density of $\frac{1}{\sqrt n}Z_n$ obtained after $n=2000$ jumps. Embedded is also a two-dimensional Gaussian probability density which has the same covariance matrix as the numerical data. Despite of the fact that the density function itself does not converge to any function, the corresponding cumulative distribution function shown in Figure~\ref{fig:5} is smooth and matches that of the corresponding Gaussian distribution. The maximum absolute difference between the numerical and analytical functions is 0.017, most of which is due to the highest peak in Figure~\ref{fig:6}.  

\section{Conclusions}
We have investigated the statistical properties of a deterministic walk in a quenched one-dimensional random environment of expanding circle maps and have analytically found the drift and variance for the resulting Gaussian probability distribution. Using numerical experiments we have been able to verify our analytical predictions. We have further studied a two-dimensional model similar to the one-dimensional system where hyperbolic toral automorphisms take the place of the circle maps. Again the probability distribution turns out to be Gaussian with certain linear drift and variance. The key feature and complicating factor in both the one- and two-dimensional cases is the \emph{fixed} random environment. A direct consequence of this is that, even after the proper scaling, the probability density does not converge to any function---a result which persists both in our one- and two-dimensional models. The implementation of recurrence to our model will be left for future work.

%%%%%%%%%%       References      %%%%%%%%%

\end{document}